\documentclass{amsart}

 \def\C{{\bf C}}  \def\N{{\bf N}} 
 \def\Q{{\bf Q}} \def\R{{\bf R}} \def\Z{{\bf Z}}

\def\End{{\mathrm{End}}}
\def\EndQ{{\mathrm{End_{\Q}}}}

\def\M{{\mathrm{M}}}
\def\NS{{\mathrm{NS}}}

\newtheorem*{Proposition}{Proposition}

\newtheorem{Corollary}{Corollary}

\theoremstyle{definition}

\newtheorem{Example}{Example}

\theoremstyle{remark}

\begin{document}

\title{Two-dimensional complex tori with multiplication by
$\sqrt{d}$}
\author{Wolfgang M. Ruppert}
\address{Mathematisches Institut, Universit\"at
Erlangen-N\"urnberg, Bismarckstra{\ss}e $1\frac{1}{2}$, D-91054
Erlangen, Germany}
\email{ruppert@mi.uni-erlangen.de}
\thanks{Version of July 7, 1998}

\begin{abstract} We give an elementary argument for the well
known fact that the endomorphism algebra $\End(A)\otimes\Q$ of a 
simple complex abelian surface $A$ can neither be an imaginary
quadratic field nor a definite quaternion algebra. Another
consequence of our argument is that a two-dimensional complex
torus $T$ with $\Q(\sqrt{d})\hookrightarrow \EndQ(T)$ where
$\Q(\sqrt{d})$ is real quadratic, is algebraic. 
\end{abstract}

\maketitle



Let $A$ be a simple abelian surface defined over $\C$. Then the
endomorphism algebra $\EndQ(A)=\End(A)\otimes_{\Z}\Q$ is one of
the following:
\begin{itemize}
\item $\Q$, 
\item a real quadratic field,
\item a CM-field of degree $4$ over $\Q$, 
\item an indefinite quaternion algebra.  
\end{itemize}
The astonishing fact is that $\EndQ(A)$ can neither be an imaginary
quadratic field nor a definite quaternion algebra, cases which
are not excluded in the usual lists for endomorphism algebras, 
cf. \cite[p.202]{M} or \cite[p.141]{LB}. In
particular if $\Q(\sqrt{d})\subseteq\EndQ(A)$ for some $d<0$
then $\EndQ(A)$ is an indefinite quaternion algebra. This fact
itself is well known (cf. \cite[Exercises (1) and (4), p.286]{LB},
\cite{S2} or \cite{OZ}). The
aim of this short note is to give an elementary proof for it. 

\bigskip

Let $\Lambda$ be a lattice in $\C^2$ and $A=\C^2/\Lambda$ the
corresponding torus. We say that $A$ (or $\Lambda$) admits
multiplication by $\sqrt{d}$ for some
$d\in\Z\setminus\{n^2:n\in\Z\}$ if there is an injective ring
homomorphism 
$$\Z[\sqrt{d}]\hookrightarrow \End(A)=\End(\Lambda),$$
i.e. if there is a matrix $D\in\M_2(\C)$ with $D^2=d$ and
$D\Lambda\subseteq \Lambda$. 

After coordinate change in $\C^2$ there are two possibilities:
$$D\sim\pm\left(\begin{array}{cc}\sqrt{d}&0\\0&\sqrt{d}
\end{array}\right)\quad\mbox{ or }\quad 
D\sim\pm\left(\begin{array}{cc}\sqrt{d}&0\\0&-\sqrt{d}
\end{array}\right).$$
The first case where $D$ is a scalar matrix, is easy: 
\begin{itemize}
\item If $d>0$ then $D$ does not operate discretely, so this
case is not possible. 
\item If $d<0$ then $\C^2/\Lambda$ is isogenous to a product of
elliptic curves with multiplication by $\sqrt{d}$ and therefore
$\EndQ(A)=\M_2(\Q(\sqrt{d}))$. 
\end{itemize}
Therefore we restrict our attention to the second case where $D$
is a nonscalar matrix. 

The N\'eron-Severi group $\NS(A)$ can be defined as the set of
hermitian forms $H(x,y)$ on $\C^2$ such that 
$\mathrm{Im}H(\Lambda,\Lambda)\subseteq \Z$ (cf. \cite[p.29]{LB}). 

Now something interesting happens:

\begin{Proposition} Let $\Lambda$ be a lattice in $\C^2$,
$D\in\M_2(\C)$ a nonscalar matrix with $D\Lambda\subseteq
\Lambda$, $D^2=d$ and
$d\in\Z\setminus\{n^2:n\in\Z\}$. Define 
$$N_D=\{H(x,y)\in \NS(\C^2/\Lambda)\mbox{ such that }H(x,Dy)\in 
\NS(\C^2/\Lambda)\} \subseteq\NS(\C^2/\Lambda).$$ 
Then $N_D\simeq\Z^2$ and 
\begin{itemize} 
\item if $d>0$ then $N_D$ contains positive definite hermitian
forms,
\item if $d<0$ then $N_D$ contains no positive definite
hermitian form. 
\end{itemize}
\end{Proposition}

The proof of the proposition is elementary and will be given in
the next section. There are some immediate consequences of the
proposition: 

\begin{Corollary} \label{real} If $\C^2/\Lambda$ admits
multiplication by $\sqrt{d}$ with $d>0$ then $\C^2/\Lambda$ is
algebraic, i.e. an abelian surface. 
\end{Corollary}

This is clear as $N_D\subseteq\NS(A)$ contains positive definite 
hermitian forms. The fact itself is also well known, cf.
\cite{S1}. 

\begin{Corollary} If $A=\C^2/\Lambda$ is an abelian surface and
admits multiplication by $\sqrt{d}$ with $d<0$ then $\NS(A)$ has
rank $\ge 3$.
\end{Corollary}

\begin{proof} If $D$ is a scalar matrix then $A$ is isogenous to
a product of elliptic curves and $\NS(A)\simeq\Z^4$. So we can
restrict us to the case that $D$ is a nonscalar matrix. 
As $A$ is abelian there is a positive definite
hermitian form $H_0\in\NS(A)$. By our proposition $H_0\not\in
N_D$. Therefore $\Z H_0\oplus N_D\subseteq\NS(A)$ and the claim
follows. 
\end{proof}

\begin{Corollary} If $A$ is an abelian surface which admits
multiplication by $\sqrt{d}$ for some $d<0$ then $A$ admits also
multiplication by $\sqrt{d'}$ for some $d'>0$ (not a square). In
particular $\EndQ(A)$ is neither an imaginary quadratic field
nor a definite quaternion algebra.
\end{Corollary}

\begin{proof} (For the following general statements see \cite[chapter
5]{LB}.) Fix a positive definite hermitian form
$H_0\in\NS(A)$. Then $\phi_{H_0}:A\to \hat{A}$ is an isogeny of
$A$ to the dual abelian variety $\hat{A}$. By 
$$\alpha\mapsto \alpha'=\phi_{H_0}^{-1}\hat{\alpha}\phi_{H_0}$$
one gets the so called Rosati involution on $\EndQ(A)$. The map
$$\NS_{\Q}(A)\to \End_{\Q}^s(A),\quad H\mapsto \phi_{H_0}^{-1}\phi_H$$
is an isomorphism (of $\Q$-vector spaces) where $s$ means
symmetric with respect to the Rosati involution. Therefore
$\dim_{\Q}\End_{\Q}^s(A)\ge 3$. The elements of
$\End_{\Q}^s(A)$ satisfy quadratic equations over $\Q$ with
positive discriminants. Therefore it is easy to find a
$d'\in\Z\setminus\{n^2:n\in\Z\}$ with $d'>0$ and
$\Q(\sqrt{d'})\hookrightarrow\EndQ(A)$. As a definite
quaternion algebra contains only imaginary quadratic subfields 
the rest is clear. 
\end{proof}

The following two examples show explicitly 
that imaginary quadratic fields
and definite quaternion algebras can be realized as endomorphism
algebras of two-dimensional complex tori. 

\begin{Example} Suppose we have $m\in\N$ and $r\in\R$ such that
$r^2,r\sqrt{m},1$ are linearly independent over $\Q$. Define 
$$\Lambda=\left(\begin{array}{cccc}
1&1+ri&\sqrt{-m}&\sqrt{-m}(1+ri)\\
1&ri&-\sqrt{-m}&-\sqrt{-m}ri
\end{array}\right)\cdot \Z^4.$$
Then it is not difficult to calculate
$\End(\Lambda)$: 
$$\End(\Lambda)=\{\left(\begin{array}{cc}
n_1+n_2\sqrt{-m}&0\\ 0&n_1-n_2\sqrt{-m}\end{array}\right):
n_1,n_2\in\Z\}\simeq\Z[\sqrt{-m}].$$
Therefore $\EndQ(\C^2/\Lambda)$ is the imaginary quadratic field
$\Q(\sqrt{-m})$. 
\end{Example}

\begin{Example} Let 
$$\Lambda=\left(\begin{array}{cccc}
1&1+\sqrt{-n}&\sqrt{-m}&\sqrt{-m}(1+\sqrt{-n})\\
1&\sqrt{-n}&-\sqrt{-m}&-\sqrt{-m}\sqrt{-n}
\end{array}\right)\cdot \Z^4$$
with $m,n\in\Z$ and $m,n\ge 1$ such that $mn$ is not a square in
$\Z$. Explicit calculation shows that 
$$\End(\C^2/\Lambda)=\Z+\Z I+\Z J+ \Z K$$
with 
$$I=\left(\begin{array}{cc}\sqrt{-m}&0\\0&-\sqrt{-m}
\end{array}\right), \quad 
J=\left(\begin{array}{cc}0&1+2\sqrt{-n}\\-1+2\sqrt{-n}&0\end{array}
\right), \quad K=IJ$$
and 
$$I^2=-m,\quad J^2=-1-4n,\quad IJ=-JI.$$
So we see that $\EndQ(\C^2/\Lambda)$ is a definite quaternion algebra. 
\end{Example}

\section*{Proof of the proposition}

Let $\Lambda$ be a lattice in $\C^2$ and let $D\in\M_2(\C)$ be a
nonscalar matrix such that $D^2=d$ for some 
$d\in\Z\setminus\{n^2:n\in\Z\}$ and $D\Lambda\subseteq\Lambda$.
After coordinate change in $\C^2$ we can assume that 
$$D=\pm\left(\begin{array}{cc}\sqrt{d}&0\\0&-\sqrt{d}\end{array}
\right).$$

\bigskip

A hermitian form $H$ on $\C^2$ is given by a hermitian matrix
$M\in \M_2(\C)$ such that 
$$\overline{M}^t=M\quad\mbox{ and }\quad
H(x,y)=x^tM\overline{y}.$$

Let $H(x,y)=x^tM\overline{y}$ be a hermitian form. Then 
$H(x,Dy)=x^tM\overline{D}\cdot \overline{y}$ is hermitian iff 
$\varepsilon MD=M\overline{D}=\overline{(M\overline{D})}^t=DM$ 
where $\varepsilon$ is the sign of $d$. This implies that $M$ is
a matrix 
\begin{eqnarray*}
M_{a,b}&=&\left(\begin{array}{cc}a&0\\0&b\end{array}\right)\mbox{ with
}a,b\in\R\mbox{ if }d>0,\\
M_{a,b}&=&\left(\begin{array}{cc}0&a+ib\\a-ib&0\end{array}\right)
\mbox{ with }a,b\in\R\mbox{ if }d<0.
\end{eqnarray*}
Note that for $d<0$ none of the matrices $M_{a,b}$ is positive
definite whereas in case $d>0$ there are positive definite
matrices $M_{a,b}$. 

\bigskip

$\Lambda$ is a $\Z[\sqrt{d}]$-module, so $\Lambda_{\Q}=
\Lambda\otimes_{\Z}\Q$ is a $\Q(\sqrt{d})$-vector space. Let
$e_1,e_2\in\Lambda$ be a $\Q(\sqrt{d})$-basis of $\Lambda_{\Q}$.
Then 
$$\Lambda_{\Q}=\Q e_1+\Q e_2+\Q De_1+\Q De_2$$
which implies that 
$$\tilde{\Lambda}=\Z e_1+\Z e_2+\Z De_1+\Z De_2\subseteq\Lambda$$
is of finite index in $\Lambda$. 

\bigskip

Write $H_{a,b}(x,y)=x^tM_{a,b}\overline{y}$ and $E_{a,b}(x,y)=
\mathrm{Im}H_{a,b}(x,y)$. Then $H_{a,b}(x,Dy)$ is also hermitian
and $E_{a,b}(x,Dy)$ is alternating. If we write 
$E_{a,b}(e_1,e_2)=u$, $E_{a,b}(e_1,De_2)=v$ then it is easy to see that 
\begin{eqnarray*}
E_{a,b}(e_1,e_2)&=&u\\
E_{a,b}(e_1,De_1)&=&0\\
E_{a,b}(e_1,De_2)&=&v\\
E_{a,b}(e_2,De_1)&=&-v\\
E_{a,b}(e_2,De_2)&=&0\\
E_{a,b}(De_1,De_2)&=&du
\end{eqnarray*}
This implies at once: If $u,v\in\Q$ then $H_{a,b}\in
\NS_{\Q}(\C^2/\Lambda)$. Define 
$$\lambda:\R^2\to\R^2,\quad (a,b)\mapsto 
(E_{a,b}(e_1,e_2),E_{a,b}(e_1,De_2)).$$
$\lambda$ is $\R$-linear and the above formulas show that $\lambda$
is injective and therefore bijective. Therefore we get an
injection
$$\Q^2\hookrightarrow \NS_{\Q}(\C^2/\Lambda),\quad z\mapsto
H_{\lambda^{-1}z}.$$
This proves the proposition.

\end{document}